\documentclass[preprint,12pt]{elsarticle}

\usepackage{lipsum}
\makeatletter 
\def\ps@pprintTitle{%
 \let\@oddhead\@empty 
 \let\@evenhead\@empty
 \def\@oddfoot{}%
 \let\@evenfoot\@oddfoot}
\makeatother
\usepackage{graphicx}
\usepackage{amssymb}
\usepackage{amsthm}
\usepackage{url}
\usepackage{lineno}
\usepackage{amssymb}
\usepackage{amsmath}
\usepackage{makeidx}
\usepackage{amscd}
\usepackage{fancyhdr}
\usepackage{yfonts}
\usepackage{graphicx}
\usepackage{bbm}
\usepackage{enumitem}
\usepackage{amsthm}
\usepackage{listings}
\usepackage{fancyvrb}
\usepackage{url}
\usepackage{mathrsfs}  
\usepackage{color}
\usepackage[all,cmtip]{xy} 

\journal{:\space To be chosen}

\newtheorem{theo.}{Theorem}[section]
\newtheorem{lem.}[theo.]{Lemma}
\newtheorem{cor.}[theo.]{Corollary}
\newtheorem{prop.}[theo.]{Proposition}
\newtheorem{def.}[theo.]{Definition}
\newtheorem{ex.}[theo.]{Example}
\newtheorem{not.}[theo.]{Notation}
\newtheorem{rem.}[theo.]{Remark}

\newcommand{\F}{\ensuremath{\mathbb F}}
\newcommand{\Z}{\ensuremath{\mathbb Z}}

\newcommand{\HC}{\ensuremath{\mathcal H}}
\newcommand{\JacH}{\ensuremath{\mathcal{J}}}

\newcommand{\tw}{\texttt{\tiny{TW}}}

\newcommand{\ab}{\mathcal{A}}

\begin{document}
\begin{frontmatter}

\title{A proof of the Hasse-Weil inequality for genus $2$ \textit{\`a la} Manin}
\author{Eduardo Ru\'{i}z Duarte\footnotemark[1] and Jaap Top}

\address{Bernoulli Institute, University of Groningen, The Netherlands}

\begin{abstract}
We prove the Hasse-Weil inequality for genus $2$ curves 
given by an equation of the form $y^2=f(x)$ with $f$
a polynomial of degree $5$,
using arguments that mimic the elementary proof of the genus $1$ case obtained by 
Yu.~I.~Manin in~1956. 
\end{abstract}

\begin{keyword}
Hasse-Weil inequality \sep Hyperelliptic curve \sep Genus two curve

\end{keyword}
\end{frontmatter}

\footnotetext[1]{{This author was supported by CONACyT M\'{e}xico through the agreement CVU-440153}}
\section{Manin's proof of the Hasse inequality for genus one}
\label{S:1}
\newcommand{\ja}{\mathcal{J}}

\newcommand{\hc}{\mathcal{H}}

Recall the following theorem:
\begin{theo.}[Hasse-Weil]
Let $C$ be an algebraic curve of genus $g$ over $\F_q$ then
\[
|\#C(\F_q)-(q+1)|\leq 2g\sqrt{q}.
\]
\end{theo.}
The Hasse-Weil inequality for an elliptic curve $E/\F_q$ (so the case of genus one, due to Hasse
in a series of papers \cite{Hasse1933}, \cite{Hasse1934}, \cite{hasse1936}) 
was obtained in an elementary way by Manin for $\text{char}\;\F_q\neq 2,3$ in \cite{manin1956oncubic} (see \cite{manin1960oncubic} for an English translation). Adjustments of these elementary arguments to the remaining 
genus one cases are presented in \cite{chahal2014supplement} and \cite{ACOFF}. Our goal is to extend these ideas to the case of genus $2$ curves. To facilitate this, we very briefly summarize in this first section Manin's argument in the genus one situation.
Throughout the paper, we restrict to finite fields of odd characteristic.\\

Let $E/\F_q$ be an elliptic curve given 
by an equation $y^2=f(x)$ where $f$ is a
polynomial of degree $3$. Consider $\phi,[n]\in\text{\rm End}_{\F_q}(E)$ where $\phi$ is the $q$-{th} Frobenius and $[n]$
the multiplication by $n$. Consider 
$\psi_n:=\phi+[n]\in \text{\rm End}_{\mathbb{F}_q}(E)$. If $\psi_n$ is non-trivial then it is of the form $(x,y)\mapsto \big(\frac{u_{1,n}(x)}{u_{2,n}(x)}, y\frac{v_1(x)}{v_2(x)}\big)$, with $u_{1,n},u_{2,n},v_1,v_2\in\F_q[x]$ such that $\text{gcd}(u_{1,n},u_{2,n})=1$ (see \cite[\S 2.9]{wash}). The Hasse-Weil inequality for $E$ follows from the claim that $d_n:=\deg(\psi_n)$
satisfies
\begin{equation}\label{hassedeg}
0\leq d_n=
\left\{\begin{array}{lr}\deg(u_{1,n})&(\mbox{if}\;\;\psi_n\neq 0)\\0&(\mbox{if}\;\;\psi_n=0)\end{array}\right.
=n^2 + (q+1-\#E(\F_q))n + q.
\end{equation}
Here for $\psi_n\neq 0$ by definition $\deg(\psi_n)=[\F_q(E):\psi_n^*\F_q(E)]$, and $\deg(u_{1,n})$ is the degree of the polynomial $u_{1,n}\in\F_q[x]$ (and $\deg(0):=0$). \\ \\
The leftmost equality in (\ref{hassedeg}) 
for $\psi_n\neq 0$ follows from the elementary observation
$\deg(\psi_n)=[\F_q(x):\F_q(\frac{u_1,n(x)}{u_2,n(x)})]=\max\{\deg (u_{1,n}(x)),\deg(u_{2,n}(x))\}$  (see \cite[\S 2.9]{wash}) or \cite[Lemma~6.2]{Soeten}), together with
$\deg(u_{1,n}(x))>\deg(u_{2,n}(x))$. The latter follows from $\psi_n(\infty)=\infty$, implying $v_\infty(\tfrac{u_{1,n}(x)}{u_{2,n}(x)})< 0$. 
The rightmost equality in (\ref{hassedeg}) is shown by induction on $n$ using the \textit{basic identity} $d_{n-1} + d_{n+1} = 2d_n + 2$ (see \cite[Lemma 8.5]{chahal1988equations}). 

Finally the non-negativity of $d_n=n^2+(q+1-\#E(\F_q))n+q$ yields that the discriminant of this quadratic polynomial in $n$ is non-positive, implying the Hasse inequality. \\

In order to extend these ideas to genus $2$ curves, we introduce an analogous $\delta_n$ which also satisfies a \textit{basic identity}, namely $\delta_{n-1}+\delta_{n+1}=2\delta_n+4$,
in the genus $2$ case.  

\section{An analogous $\delta_n$ for genus $2$}
Let $k:=\F_q$ be a finite field of odd cardinality $q$, and let $\HC$ be a hyperelliptic curve of genus $2$ over $k$. Throughout, we assume that $\HC$ is given by the equation $Y^2=X^5 +a_4X^4+a_3X^3+a_2X^2+a_1X+a_0$. By $\JacH$ we denote the Jacobian variety associated to $\HC$. The points of $\JacH$ correspond to divisor classes $[D]\in\text{Pic}^0(\HC)$. 

Fix the point $\infty\in\HC$ and consider the Abel-Jacobi map $\iota\in \text{Mor}_k(\HC,\JacH)$ given by $P\mapsto [P-\infty]$. We have that $\Theta:=\text{Im}\;\iota$ is the \textit{theta divisor} of $\JacH$ and $\Theta\cong\HC$. Consider the $q$-th power 
Frobenius map $\phi\in\text{End}_{k}(\JacH)$ and the morphism  $\Phi:=\phi\circ\iota\in\text{Mor}_{k}(\HC,\JacH)$. Since $\text{Mor}_k(\HC,\JacH)\cong \JacH(k(\HC))$ is an Abelian group, we define $\Psi_n:=\Phi+n\cdot\iota\in\text{Mor}_k(\HC,\JacH)$ and $\Theta_n:=\text{Im}\;\Psi_n\subset\JacH $. \\

Similar to $d_n=\deg(\psi_n)$ in (\ref{hassedeg}), we will define the
`complexity' (height) of $\Psi_n$ and denote this by $\delta_n$. 

Consider the generic point of $\JacH$ given by $\mathfrak{g}:=[(x_1,y_1)+(x_2,y_2)-2\infty]$ and the function field of $\JacH$ denoted by $k(\JacH)\cong k(x_1+x_2,x_1x_2,\tfrac{y_1-y_2}{x_1-x_2},\tfrac{x_2y_1-x_1y_2}{x_1-x_2})$. The Riemann-Roch space $\mathcal{L}(2\Theta)\subset k(\JacH)$ has dimension $4$; a basis is 
given by $\{\kappa_1,\kappa_2,\kappa_3,\kappa_4\}$ with $\kappa_1:=1, \kappa_2:=x_1+x_2, \kappa_3:=x_1x_2, \kappa_4:=\tfrac{F_0(x_1+x_2,x_1x_2)-2y_1y_2}{(x_1-x_2)^2}$ where $F_0(A,B):= 2a_0 + a_1A+2a_2B+a_3AB+2a_4B^2 + AB^2$ (see \cite[Chapter 2]{PROLEGOMENA} and \cite[Page 5]{FLYNNJAC}). This basis 
is used to define a singular surface $\mathcal{K}\!\subset\!\mathbb{P}^3$ birational to the Kummer surface associated to~$\JacH$, as the closure of the image of
\begin{align*}
\kappa:\JacH\setminus \Theta &\to \mathbb{P}^3,\\
D&\mapsto [\kappa_1(D):\kappa_2(D):\kappa_3(D):\kappa_4(D)].
\end{align*}
As a \textbf{remark}, let $[-1]\in\text{Aut}(\JacH)$ be the involution on $\JacH$ given by $[-1]\mathfrak{g}=[(x_1,-y_1)+(x_2,-y_2)-2\infty]=-\mathfrak{g}$. Note that $[-1]^*:k(\JacH)\to k(\JacH)$ is the trivial map on $\mathcal{L}(2\Theta)\subset k(\JacH)$. Particularly for all $D\in\ja\setminus\Theta$ we have that $\kappa_i(D)=\kappa_i(-D)$.

Using the previous remark, suppose that $\Theta_n:=\text{Im}\;\Psi_n\not\subset \Theta$ and $\Theta_n$ is not a zero of $\kappa_4$. Let $(x,y)\in\HC$ be generic, then we have that $\psi_n(x,y):=(\kappa_4\circ\Psi_n)(x,y)=(\kappa_4\circ\Psi_n)(x,-y)$. Hence $\psi_n(x,y)\in k(x)$, that is, $\psi_n(x,y)=:\tfrac{\mu_{1,n}(x)}{\mu_{2,n}(x)}$ is a rational function in one variable $x$, which is the one that we will use to define $\delta_n$.  \\

Our geometric situation is described in the following diagram:

\begin{equation}\label{g2diag2}
\xymatrix{
 &\JacH \ar[dr]^{\phi+[n]}\\
\HC \ar[ur]^{\iota}  \ar[ddrr]_{\psi_n} \ar[rr]^{\Psi_n} &&  \JacH \ar@{.>}[d]^{\kappa}_{2:1} \ar@{.>}@/^2pc/[dd]^{\kappa_4}&   \\
                     &&\mathcal{K} \ar@{.>}[d]^{\pi}\\
                     &&\mathbb{P}^1
} 
\end{equation}
Here $\pi$ is a projection and $\kappa_4=\pi\circ\kappa$. Since $\Theta_n$ is not a zero or a pole of $\kappa_4$, we define $\delta_n:=\deg\mu_{1,n}$.\\

There are two situations left to define $\delta_n$ for every $n\in\Z$. The first is when $\Psi_n$ is constant (hence equal to the zero map). In this case $\Theta_n\subset\ja$ is a point and we define $\delta_n:=0$. The second is when $\{0\}\neq \Theta_n$ is a curve which is a zero or a pole of $\kappa_4$, that is $\Theta_n\in \text{Supp div}(\kappa_4)$. In the following section (Formula~(\ref{deltan})) we will define $\delta_n$ for this special situation. We show that if $\Psi_n$ is non-constant but  $\kappa_4(\Psi_n(x,y))=c$ is constant then $c$ can only be $0$ or $\infty$, that is, the curve $\Theta_n=\text{Im}\Psi_n$ is a zero or a pole of $\kappa_4$ respectively (see Lemma \ref{k4constant}).  \\

Further, if $\Theta_n$ is not a zero or a pole of $\kappa_4$, we show 
in the next section that $\tfrac{\deg \psi_n}{2}=\max\{\deg\mu_{1,n},\deg\mu_{2,n}\}$.
In the case $\Theta_n$ is a zero or a pole of $\kappa_4$, a similar equality will be shown taking a translation of $\Theta_n$ by a $2$-torsion point of $\ja$ in order to avoid the pole and zero divisor of $\ja$.\\

Finally, we show the \textit{basic identity} $\delta_{n-1}+\delta_{n+1}=2\delta_n+4$. The same strategy also employed by Manin for genus one
will then lead to a proof of the Hasse-Weil inequality in this case.

\section{Proof of the Hasse-Weil inequality for genus $2$}
We use the notations $\HC/\F_q$, $\Theta,\Theta_n$ etc.\ introduced in the previous section.  
\begin{lem.}\label{maxpsi}
Let $(x,y)\in\HC/\F_q$ be generic. Suppose that $\Theta_n\subset\ja$ is a curve that is not a zero nor a pole of $\kappa_4$. Then $\psi_n(x,y)=\tfrac{\mu_{1,n}(x)}{\mu_{2,n}(x)}$ for some
coprime polynomials $\mu_{1,n},\mu_{2,n}$. \\Moreover $\tfrac{\deg\psi_n}{2}=\max\{\deg\mu_{1,n}(x),\deg\mu_{2,n}(x)\}$.
\end{lem.}
\begin{proof}
Since $\Theta_n$ is not a zero or a pole of $\kappa_4$ and $\kappa_4\in\mathcal{L}(2\Theta)$, the function $\psi_n(x,y)=\kappa_4(\Psi_n(x,y))$ is defined and non-zero. In the previous section we saw that $\psi_n(x,y)=\psi_n(x,-y)$, and hence $\psi_n(x,y)\in \F_q(x)$. This shows
the existence of the coprime $\mu_{1,n},\mu_{2,n}\in \F_q[x]$.

Now,
\[
\deg(\psi_n)=[\F_q(x,y):\F_q(\tfrac{\mu_{1,n}(x)}{\mu_{2,n}(x)})]=
[\F_q(x,y):\F_q(x)]\cdot [\F_q(x):\F_q(\tfrac{\mu_{1,n}(x)}{\mu_{2,n}(x)})]
\]
and the lemma follows.
\end{proof}

It can happen that $\Psi_n$ is the zero map, which implies $\Theta_n=\text{Im}\Psi_n$ is a point. For example consider the hyperelliptic curve $Y^2 = X^5 + 5X$ over $\F_{49}$. An explicit MAGMA computation shows that $\psi_7:=\Phi+7\iota\in\text{\rm Mor}_{\F_{49}}(\HC,\JacH)\cong\ja(\F_{49}(\HC))$ is the zero map: 

\begin{Verbatim}[fontsize=\scriptsize]
> p := 7; q := p^2; F := FiniteField(q);
> P<x> := PolynomialRing(F);
> H := HyperellipticCurve(x^5+5*x); 
> FH<X,Y> := FunctionField(H); HE := BaseExtend(H,FH);
> JE := Jacobian(HE); M<t> := PolynomialRing(FH);      
> Phi := JE![t-X^q, Y^q];
> GPt := JE![t-X,   Y];
> -7*GPt;
(x + 6*X^49, (X^120 + X^116 + 5*X^112 + 6*X^108 + X^92 + X^88 + 5*X^84 + 6*X^80 + 5*X^64 + 
5*X^60 + 4*X^56 + 2*X^52 + 6*X^36 + 6*X^32 + 2*X^28 + X^24)*Y, 1)
> Phi;
(x + 6*X^49, (X^120 + X^116 + 5*X^112 + 6*X^108 + X^92 + X^88 + 5*X^84 + 6*X^80 + 5*X^64 + 
5*X^60 + 4*X^56 + 2*X^52 + 6*X^36 + 6*X^32 + 2*X^28 + X^24)*Y, 1)
> Phi+7*GPt;
(1, 0, 0)
\end{Verbatim}
In this example $\JacH$ is isogenous to the square of a supersingular elliptic curve and the ground field has $p^2=49$ elements. The characteristic polynomial of Frobenius $\phi\in\text{\rm End}_{\F_{49}}(\JacH)$ is given by $\chi_\phi(X):=({X}+7)^4$ which is the main reason of this behavior. \\
A general construction of  curves having  Jacobian isogenous to a square of a supersingular elliptic curve was achieved by Moret-Bailly in \cite{moretcurves}.\\

The following proposition and lemma isolates a special case for our final proof of the Hasse-Weil inequality for genus $2$. Note that the example discussed above illustrates
this special case.
\begin{prop.}\label{exceptional}
Let $\HC/\F_q$ be a hyperelliptic curve of genus $2$,
given by an equation $y^2=f(x)$ with $f$ of degree $5$.
Let $\JacH$ be the Jacobian of $\HC$ and
$\iota\colon \HC\to\JacH$ the map $P\mapsto [P-\infty]$. Suppose that there is an $n\in\Z$ such that $\Psi_n=(\phi+[n])\circ\iota\in\text{\rm Mor}_{\F_q}(\HC,\JacH)$ is constant. Then $q$ is a perfect square and 
$\#\HC(\F_q)=q+1+4n=q+1\pm 4\sqrt{q}$.
\end{prop.}
\begin{proof}
First, we show that if $\Psi_n=(\phi+[n])\circ\iota$ is constant, then $\phi=-[n]$.\\
We have that $\Psi_n=(\phi+[n])\circ\iota$ is constant and $0\in\text{Im}\Psi_n$, hence $\Psi_n=0$; this is equivalent to $(\phi+[n])(\Theta)=0$ since  $\Theta=\iota(\HC)$. Moreover, $\Theta$ generates $\JacH$, that is $\JacH=\{D_1+D_2: D_1,D_2\in\Theta\}$.
So if any $\varphi\in\mbox{End}(\JacH)$ is zero on
$\Theta$ then it is the zero map.
Hence $\phi=-[n]\in\text{End}(\JacH)$. \\
Note that $\phi=-[n]$ implies
$q^2=\deg(\phi)=\deg([-n])=n^4$. Hence  $q=n^2$ is a perfect square and $n=\pm \sqrt{q}$. \\ \\
Now we proceed to count $\#\HC(\F_q)$. Using that $\phi=-[n]$ we have that:
\begin{equation}\label{ptsexc}
\#\JacH(\F_q)=\#\text{Ker}(\phi-[1])=\#\text{Ker}(-[n+1])=(n+1)^{4}.
\end{equation}
(Here we used that $n+1$ is not a multiple of $\mbox{char}(\F_q)$).
Moreover, an easy counting argument (see \cite[Chapter 8, \S 2]{PROLEGOMENA}) shows:
\begin{equation}\label{countjac}
\#\JacH(\F_q)=\displaystyle\frac{\left(\#\HC(\F_q)\right)^2 + \#\HC(\F_{q^2})}{2}-q.
\end{equation}
Consider the quadratic twist of $\HC$ denoted by $\HC^{\texttt{\tiny{TW}}}$ and its Jacobian  $\JacH^{\tw}$, then:
\begin{equation}\label{ptsexctw}
\#\JacH^{\tw}(\F_q)=\#\text{Ker}(\phi+[1])=\#\text{Ker}(-[n]+1)=(n-1)^4.
\end{equation}
Similarly as in (\ref{countjac}) and using that $\#\HC(\F_q)+\#\HC^\tw(\F_q)=2q+2=2n^2+2$ and $\HC^{\tw}(\F_{q^2})\cong \HC(\F_{q^2})$, we have that:
\begin{align}\label{countjactw}
\begin{split}
\#\JacH^{\tw}(\F_q)&=\displaystyle\frac{\#\HC^{\tw}(\F_q)^2 + \#\HC(\F_{q^2})}{2}-q\\
&=\displaystyle\frac{(2n^2+2-\#\HC(\F_q))^2 + \#\HC(\F_{q^2})}{2}-q=(n-1)^4
\end{split}
\end{align}

Subtracting  (\ref{countjactw})
from (\ref{countjac}) yields:
\begin{align}\label{finalexc}
\begin{split}
\#\HC(\F_q)^2 - (2n^2+2-\#\HC(\F_q))^2&=2\big((n+1)^4 - (n-1)^4\big)\\
&=16n(n^2+1),
\end{split}
\end{align}
which can be rewritten as $\#\HC(\F_q)=n^2+4n+1=q+1\pm 4\sqrt{q}$.
\end{proof}
\begin{cor.}\label{constantmap}
Let $\HC/\F_q$ be a hyperelliptic curve of genus $2$ given by an equation $Y^2 = f(X)$ with $f$ of degree $5$. Let $\ja$ be its Jacobian and suppose that $\Psi_n=(\phi+[n])\circ\iota$ is constant. Then $\text{\rm Im}\;\Psi_{n-1}=\text{\rm Im}\;\Psi_{n+1}=\Theta\in\text{\rm Div}(\ja)$.  
\end{cor.}
\begin{proof}
By the Proposition (\ref{exceptional}) $\phi=-[n]$, hence $\Psi_{n\pm 1}=\pm [1]\circ\iota=\pm \iota\in\text{Mor}_{\F_q}(\HC,\JacH)$. Further, $\text{Im}\;\iota=\text{Im}-\!\iota=\Theta$ since $\Theta\cong\HC$ (here we used that $\Theta$ is symmetric with respect of the hyperelliptic involution under $\iota$).   
\end{proof}
\noindent
\begin{lem.}\label{invfi}
Let $\HC/\F_q$ be a hyperelliptic curve of genus $2$ given by $Y^2=f(X)$ with $\deg f(X)=5$ and let $\JacH$ be the
Jacobian of $\HC$. Let $(x,y)\in\HC$ be generic, then $-\Psi_n(x,y)=\Psi_n(x,-y)$.
\end{lem.}
\begin{proof}
Using $[(a,b)-\infty]=[\infty-(a,-b)]$ for any $(a,b)\in\HC$ one finds
\begin{align*}
-\Psi_n(x,y)&=-\phi([(x,y)-\infty])-n[(x,y)-\infty]\\
&=[\infty-(x^q,y^q)]+ n[\infty-(x,y)]\\
&=[(x^q,-y^q)-\infty]+n[(x,-y)-\infty]=\Psi_n(x,-y).
\end{align*}
\end{proof}
Now we  calculate some values of $\delta_n$.
\begin{prop.}\label{k41-1}
Let $\HC/\F_q:\;Y^2=X^5+a_4X^4 + a_3X^3+a_2X^2+a_1X+a_0=:f(X)$ be a genus $2$ curve ($q$ odd).  Then $\delta_{-1}=\#\HC(\F_q)+q+1$.
\end{prop.}
\begin{proof}
For $(x,y)\in\HC$ generic, $\Psi_{-1}(x,y)=[(x^q,yf(x)^{\frac{q-1}{2}})+(x,-y)-2\infty]$.
Lemma \ref{maxpsi} shows that $\delta_{-1}$ equals the degree of $\psi_{-1}(x,y)=\kappa_4(\Psi_{-1}(x,y))\in \F_q(x)$. This degree is the maximum of the polynomial degrees of the numerator and the denominator, assuming these are coprime. Here
\begin{equation}\label{k4fi}
\psi_{-1}(x,y)=
\tfrac{x^{3q+2}+x^{2q+3}+2a_4x^{2q+2}+a_3(x^{2q+1}+x^{q+2})+2a_2x^{q+1}+a_1(x^q+x)+2a_0+
2f(x)^{\frac{q+1}{2}} }{(x^q-x)^2}.
\end{equation}
Let $\nu(x)$ and $\eta(x)$ be respectively the numerator and denominator of $(\ref{k4fi})$ before cancellation of common factors.  Note that $\deg(\eta)=2q$ and that every $\alpha\in\F_q$ is a double root of $\eta(x)$. Furthermore, $\deg(\nu(x))=3q+2>2q$, hence $\delta_{-1}=3q+2-\deg(\gcd(\nu(x),\eta(x)))$.\\
Since $\psi_{-1}(x,y)\in\F_q(\HC)$, the common factors $(x-\alpha)$ of $\nu$ and $\eta$ occur at the points $(\alpha,\beta)\in\HC$ such that $\alpha\in\F_q$ and $\beta\in\F_q^*$ or $\beta\in\F_{q^2}^*\setminus\F_{q}^*$ or $\beta=0$. Hence, we have three possibilities for cancellations:    \\ \\
\textbf{Case $\beta\in\F_{q}^*$:}\\
In this case $(\alpha,\beta)\in\HC(\F_q)$ and therefore $f(\alpha)$ is a square in $\F_q^*$. Hence $f(\alpha)^{\frac{q-1}{2}}=1$. Moreover, $\alpha^q=\alpha$. Using this, the last term of $\nu(\alpha)$ is $2f(\alpha)^{\frac{q+1}{2}}=2f(\alpha)f(\alpha)^{\frac{q-1}{2}}=2f(\alpha)$ and
\[\nu(\alpha)
=4f(\alpha).
\]
Since $\beta\neq 0$ there is no cancellation of a factor $(x-\alpha)$ in this case.\\ \\
\textbf{Case $\beta=0$:}\\
We have that $f(\alpha)=0$ and $\alpha^q=\alpha$, so the numerator of (\ref{k4fi}) is $2f(\alpha)=0$. Therefore $\nu(x)$ and $\eta(x)$ share  the linear factor $x-\alpha$ with multiplicity one or two. 
The multiplicity in fact equals one since $\frac{d}{dx}\nu(x)\mid_\alpha=4f'(\alpha)\neq 0$ as $f(x)$ does not have repeated zeros. \\ \\
\textbf{Case $\beta\not\in\F_{q}$:}\\
In this case $f(\alpha)$ is nonzero and is not a square in $\F_q^*$. Therefore $f(\alpha)^{\frac{q-1}{2}}=-1$ by Euler's criterion. Moreover  $\alpha^q=\alpha$  and  $\nu(\alpha)$ is in this case 
\[2(\alpha^5+a_4\alpha^4+a_3\alpha^3+a_2\alpha^2+a_1\alpha+a_0)-2f(\alpha)
=0.
\]
 To find the multiplicity of $\alpha$ as a zero of $\nu(x)$, note that
\begin{align*}
\frac{d}{dx}\nu(x)\mid_\alpha&=2\alpha^{3q+1}\!+\!3\alpha^{2q+2}\!+\!4a_4\alpha^{2q+1}\!+\!a_3(\alpha^{2q}\!+\!2\alpha^{q+1})\!+\!2a_2\alpha^q\!+\!a_1\!-\!f'(\alpha)\\
&=5\alpha^4+ 4a_4\alpha^3+3a_3\alpha^2 + 2a_2\alpha+a_1-f'(\alpha)\\
&=f'(\alpha)-f'(\alpha).\\
&=0.
\end{align*} 
This tells us that  the factor $(x-\alpha)^2$ appears in $\nu$ and then it cancels with the denominator.

Combining the cases, one concludes
$\deg(\gcd(\nu(x),\eta(x)))=2q+1-\#\HC(\F_q)$ and therefore $\deg(\kappa_4(\Psi_{-1}(x,y)))=
q+1+\#\HC(\F_q)
$.
\end{proof}
\begin{prop.}\label{corok41}
With notations as in Proposition \ref{k41-1} one has
\[
\delta_1=3(q+1)-\#\HC(\F_q).
\]
\end{prop.}
\begin{proof}
Note that
\begin{equation*}\label{k4tw}
\kappa_4(\Psi_1(x,y))=\tfrac{x^{3q+2}+x^{2q+3}+2a_4x^{2q+2}+a_3(x^{2q+1}+x^{q+2})+2a_2x^{q+1}+a_1(x^q+x)+2a_0-2f(x)^{\frac{q+1}{2}} }{(x^q-x)^2}.
\end{equation*}
This expression differs from (\ref{k4fi}) only at the sign of the last term of the numerator, namely $2f(x)^{\frac{q+1}{2}}$. An analogous argument as the one given in Proposition \ref{k41-1} proves the proposition.
\end{proof}
We will use the following definition in order to interpret $\delta_n$ in the case where $\Theta_n$ is a curve that is a zero or a pole of $\kappa_4$. The case $\Theta_n=\Theta$ (which is a pole of $\kappa_4$) occurs, e.g., for $n=0$.
\begin{def.}
Let $D_1, D_2\in\text{\rm Div}(\JacH)$. 
By $D_1\bullet D_2$ we denote the intersection number of the
divisors  $D_1$ and $D_2$ on the surface $\JacH$.   
\end{def.}
\noindent
For details and properties of this see, e.g., \cite[Appendix C or Chapter V]{HARTSHORNE}. 
\begin{lem.}\label{degintersection}
Suppose that $\text{\rm Im}\Psi_n=\Theta_n\not\subset\Theta$
and that $\psi_n=\kappa_4\circ \Phi_n\circ\iota\colon \HC\to\mathbb{P}^1$ is
nonconstant, where $\Phi_n:=\phi+[n]\in\text{\rm End}(\JacH)$. Then 
\[
2\Theta\bullet\Theta_n={\deg\psi_n}=2\Phi_n^*\Theta\bullet\Theta.
\]
\end{lem.}
\begin{proof}
Let $(x,y)\in\HC$ be generic. Since $\text{Im}\Psi_n=\Theta_n\not\subset\Theta$ and since $\kappa_4\in
\F_q(\JacH)$ has divisor
$D-2\Theta$ for some effective divisor $D\in\text{Div}(\JacH)$, we have that $\psi_n(x,y)=\kappa_4(\Psi_n(x,y))\in\F_q(x)$ by Lemma \ref{maxpsi}.  Moreover by assumption this rational function is nonconstant. Therefore $\deg\psi_n=
\deg\left(({\kappa_4}_{\mid{\Theta_n}})^*\infty\right)
=2\Theta\bullet \Theta_n$ which shows the first equality. 

For the second, note that $\Phi_n^{-1}(\Theta)=\{D\in\JacH : \Phi_n(D)\in\Theta\}$ and $\Theta$ is the locus where $\kappa_4$
has a pole (in fact a double pole). 
Since $\deg(\psi_n)=\deg({\kappa_4\circ \Phi_n}_{\mid\Theta})$ we conclude $\deg\psi_n=2\Phi_n^{-1}(\Theta)\bullet\Theta$. 
 Applying \cite[Lemma~1.7.1]{intfulton}
 this equals $2\Phi_n^*\Theta\bullet\Theta$.
\end{proof}

We will deal with the cases where $\Theta_n=\Theta$ by using a linear equivalent divisor $\Theta_n'\in\text{Div}(\JacH)$ and Lemma \ref{degintersection}. \\

First we show that $\Theta\bullet\Theta=2$ using the pole divisor of $\kappa_4$. 
To achieve this we use (see \cite[Chapter V. Theorem 1.1]{HARTSHORNE}) that if $\Theta\sim\Theta'$ as divisors ($\sim$ denoting linear equivalence), then $\Theta\bullet\Theta=\Theta\bullet\Theta'$. A suitable divisor $\Theta'$  will be constructed as a symmetric translation of  $\Theta\subset\JacH$ (with respect of $[-1]\in\text{Aut}(\ja)$). 

\begin{rem.}\label{remtrans}
{\rm The following lemmas use that for divisors $D,D'$ on
$\ja$ and any point $\xi\in\ja$, denoting by $t_\xi$ the translation over $\xi$,
we have $D\bullet D'=D\bullet t_\xi^*D'$. Indeed, the fact that $D'$ and $t_\xi^*D'$
are algebraically equivalent is a special case of \cite[I \S{2} Proposition~2]{Lang},
and the fact that algebraically equivalent cycles are numerically equivalent
can be found in \cite[19.3.1]{intfulton}.
}\end{rem.}
\begin{lem.}\label{selftheta}
Let $\HC/\F_q$ be a hyperelliptic curve of genus $2$ given by $Y^2=X^5 + a_4X^4 + a_3X^3 + a_2X^2 + a_1X + a_0=:f(X)$ and consider its Jacobian $\JacH$.\\ 
Then $\Theta=\iota(\HC)\subset\JacH$ satisfies $\Theta\bullet\Theta=2$.
\end{lem.}
\begin{proof}
Let $(w,0)\in\HC(\overline{\F}_q)$ be a Weierstrass point and consider $\iota_w\in\text{Mor}(\HC,\JacH)$ given by $P\mapsto [P+(w,0)-2\infty]$. Let $(x,y)\in\HC$ be the generic point. We have that $\Theta':=\text{Im}\;\iota_w\subset\JacH$ is a translation of $\Theta$, and therefore $\Theta'\sim\Theta$ in $\text{Div}(\JacH)$. 
Since $\kappa_4\circ [-1]=\kappa_4$ and $[(x,-y)+(w,0)-2\infty]=2\infty-(x,y)-(w,0)]$
it follows that $\kappa_4(\iota_w(x,y))\in\F_q(x)$. Analogous to the proof
of Lemma~\ref{degintersection} one obtains $\Theta\bullet\Theta=\Theta'\bullet\Theta=\deg\kappa_4(\iota_w(x,y))$ where deg denotes
the degree of the given element of $\F_q(x)$ (which is half the degree of the
map $\kappa_4\circ\iota_w\colon \HC\to \mathbb{P}^1$). 
Note that 
\begin{equation}\label{selfwithkappa}
\kappa_4(\iota_w(x,y))=\tfrac{2a_0+a_1(x+w)+2a_2xw+a_3(x+w)xw+2a_4(xw)^2+(x+w)(xw)^2}{(x-w)^2}
\end{equation}
One observes that both the numerator and denominator here are divisible by $(x-w)$ 
and the derivative w.r.t.\ $x$ of the numerator, evaluated at $x=w$, equals $f'(w)\neq 0$. Hence $\deg\kappa_4(\iota_w(x,y))=2$ which proves the lemma. 
\end{proof}
\noindent
Lemma~\ref{selftheta} also follows using the \textit{adjunction formula} (\cite[Chapter V, 1.5]{HARTSHORNE}). \\

Let $\Phi_n:=\phi+[n]\in\text{End}(\ja)$. Using an analogous argument, we  calculate $(\Phi_0)_*\Theta\bullet\Theta
=\Theta\bullet\Phi_0^*\Theta$; the equality of these intersection numbers is a
consequence of the `projection formula' \cite[Proposition~2.3(c)]{intfulton}.
In fact, in the present case equality is also established by computing both
numbers directly.
\begin{lem.}\label{selfthetafrob}
Let $\HC/\F_q$ be a hyperelliptic curve of genus $2$ given by $Y^2=X^5 + a_4X^4 + a_3X^3 + a_2X^2 + a_1X + a_0=:f(X)$ and consider its Jacobian $\JacH$. With notations
as before, we have $(\Phi_0)_*\Theta\bullet\Theta= 2q=\Theta\bullet\Phi_0^*\Theta$.
\end{lem.}
\begin{proof}
As $\Phi_0\colon \ja\to\ja$ is the $q$th power Frobenius morphism,
its restriction to $\Theta$ maps $\Theta$ to itself and has degree $q$.
As a consequence $(\Phi_0)_*\Theta=q\Theta$ and therefore
Lemma~\ref{selftheta} implies the first equality.
For the second equality, let $(x,y)\in\HC$ be generic. By a similar argument the
one presented in Lemma~\ref{selftheta} we translate $\Psi_0\in\text{Mor}(\HC,\ja)$ by $\iota(w,0)\in \ja(\overline{\F}_q)$, namely $\Psi_{0,w}(x,y)=[(x^q,yf(x)^{\tfrac{q-1}{2}})+[(w,0)-2\infty]$. Then, since $\text{Im}\Psi_0=\Phi_0(\Theta)$ we define $\Phi_{0,w}:=\Phi_0+\iota(w,0)$ and we have that $\Theta\bullet\Phi_0^*\Theta=\Theta\bullet\Phi_{0,w}^*\Theta$. The latter intersection number equals the degree of 
\begin{equation}\label{selfwithkappafrob}
\kappa_4(\Psi_{0,w}(x,y))=\tfrac{2a_0+a_1(x^q+w)+2a_2x^q w+a_3(x^q+w)x^qw+2a_4(x^qw)^2+(x^q+w)(x^qw)^2}{(x^q-w)^2}
\end{equation}
(considered as a morphism $\mathbb{P}^1\to\mathbb{P}^1$). Take $v\in\overline{\F}_q$
with $v^q=w$, then the denominator in the right-hand-side of (\ref{selfwithkappafrob}) equals $(x-v)^{2q}$. The numerator equals
$(2a_0+a_1(x+v)+2a_2xv+a_3(x+v)xv+2a_4(xv)^2 + (x+v)xv)^q$. 
Evaluating the numerator at $x=v$ yields $(2f(v))^q=0$, hence the numerator is divisible by $(x - v)^q$. Since the derivative of
$2a_0+a_1(x+v)+2a_2xv+a_3(x+v)xv+2a_4(xv)^2 + (x+v)xv$ evaluated at $x=v$
equals $f'(v)\neq 0$ it follows that the rational function (\ref{selfwithkappafrob}) has degree $2q$. 
\end{proof}
\begin{lem.}\label{k4constant}
Suppose that $\Theta_n\subset\ja$ is a curve. If $\kappa_4(\Psi_n(x,y))=c$ is constant then $c\in \{0,\infty\}$. 
\end{lem.}
\begin{proof}
Let $\Theta_n=\text{Im}\Psi_n\in \text{Supp}\text{ div}(\kappa_4)$ then $\kappa_4(\Psi_n(x,y))\in\{0,\infty\}$ depending on $\Theta_n$ being  a zero or a pole of $\kappa_4$.\\
Suppose that $\kappa_4(\Psi_n(x,y))=c\in {\F}_q^*$. Since $\kappa_4\in\mathcal{L}(2\Theta)$ and $\deg \kappa_4(\Psi_n(x,y))=\Theta_n\bullet\Theta=0$ (note that here we use $c\neq 0,\infty$), this contradicts the fact that the curves $\Theta_n$ and $\Theta$ intersect in $0\in\ja$.   
With this we have that $c\in\{0,\infty\}$.     
\end{proof}
Lemma~\ref{k4constant} (rather, its proof)
provides a geometric
reason for the fact that $\kappa_4\in\mathcal{L}(2\Theta)$ cannot be constant 
$\neq 0,\infty$ when restricted to the curves $\Theta_n$:  these curves will always intersect $\Theta$ and therefore they will have a positive intersection number (degree of $\kappa_4(\Psi_n(x,y)$). 
However, if the curve $\Theta_n$ equals
$\Theta$ or is contained in the
zero-locus of $\kappa_4$, then 
of course $\kappa_4|_{\Theta_n}$ is the
constant map $\infty$ or $0$, respectively.

Note that when $\Theta_n$ has dimension zero, then
$\Theta_n=\{0\}\subset\ja$ (compare the proof of Proposition~\ref{exceptional}). 
This case is treated separately in the definition of $\delta_n$ given below.

Using Lemma~\ref{k4constant}, it follows that the only cases
where  $\delta_n$ is not defined yet, is the situation where $\Theta_n$ is a pole or a zero of $\kappa_4$. The following
property of $\kappa_4$ will be useful.\\
\begin{lem.}\label{zeroeskappa4_fix}
Let $\HC/\F_q$ be a hyperelliptic curve
given by an equation $Y^2=f(X)$ where
$f$ has degree $5$ and consider its Jacobian $\ja$.
Then the function $\kappa_4$ on $\ja$
has a divisor of zeros $D_0=\textrm{ \rm div}_0(\kappa_4)$ such that
its support consists of at most four irreducible curves.
\end{lem.}
\begin{proof}
Let $\text{div}(\kappa_4)=D_0-2\Theta$ with $D_0$ effective. We have that $D_0=\sum C_i$ where the $C_i\subset\ja$ are irreducible curves.
Then
\[
4=2\Theta\bullet\Theta=
\left(\sum C_i\right)\bullet\Theta=\sum(C_i\bullet\Theta).
\]
As $\Theta$ is ample, $C_i\bullet\Theta>0$ (see \cite[Chapter V, Theorem 1.10]{HARTSHORNE}) which implies that the support of the zero divisor of $\kappa_4$ consists of four or less irreducible curves in $\ja$.
\end{proof}
\color{black}
\begin{lem.}\label{thetanintheta}
Let $\HC/\F_q$ be a hyperelliptic curve of genus $2$ given by $Y^2=X^5 + a_4X^4 + a_3X^3 + a_2X^2 + a_1X + a_0=f(X)$ and let $\ja$ be its Jacobian. 
Assume $\text{\rm Im}\;\Psi_n=\Theta_n\subset\ja$ is a curve that is a zero or a pole of $\kappa_4$. Let $(w,0)\in \HC(\overline{\F}_{q})$ be some Weierstrass point. Define $\Phi_{n,w}:=\Phi_n+\iota(w,0)$ where $\Phi_n=\phi+[n]\in\mbox{End}(\ja)$ and consider  $\Psi_{n,w}:=\Psi_n + \iota(w,0)\in\text{\rm Mor}_{\overline{\F}_{q}}(\HC,\ja)$. Take the generic point $(x,y)\in\HC$, then $\kappa_4(\Psi_{n,w}(x,y))=:\tfrac{\mu_{1,n}^w(x)}{\mu_{2,n}^w(x)}\in\overline{\F}_{q}(x)^*$ for coprime
$\mu_{1,n}^w,\mu_{2,n}^w\in\overline{\F}_{q}[x]$ 
 and ${\Phi_{n,w}}^*\Theta\bullet\Theta=\max\{\deg\mu_{1,n}^w,\deg \mu_{2,n}^w\}$.
\end{lem.}
\begin{proof}
Suppose that $\Psi_n(\HC)=\Theta_n=\Theta$ (is a pole of $\kappa_4$), it follows that $\Theta_n^w:=\text{Im}\Psi_{n,w}\not\subset\Theta$. We want $\Theta_n^w$ to avoid the support of $\kappa_4$. We do this in order to have a well defined  degree of $\kappa_4$ restricted to the symmetric divisor 
$\Theta_n^w$ (see Remark \ref{remtrans}).

If the curve $\Theta_n^w\subset\ja$ would be a zero of $\kappa_4$, then by Lemma~\ref{zeroeskappa4_fix} there are at most four possible curves $\{C_1,C_2,C_3,C_4\}$ in the zero locus $\kappa_4$. Further, there are exactly five affine Weierstrass points in $\HC$, and then for at least one of them say $(\hat w,0)\in\HC$ we have that $\Theta_n^{\hat w }\not\in\{C_1,C_2,C_3,C_4\}$. \\
\color{black}

As a result,  $\kappa_4(\Psi_{n,{\hat w}}(x,y))$ is well defined and non-constant since $\Theta_n^{\hat w}\not\in\text{Supp}\text{ div}(\kappa_4)$ (see Lemma~\ref{k4constant}). Further, $\text{Im}\Psi_{n,\hat w}$ is symmetric with respect to $[-1]\in\text{Aut}(\ja)$ since $\Psi_n(x,y)+\iota( {\hat w},0)$ is. Hence $\kappa_4(\Psi_{n,{\hat w}}(x,y))=\kappa_4(\Psi_{n, {\hat w}}(x,-y))\in\overline \F_{q}(x)^*$. 

By Remark \ref{remtrans} and Lemma \ref{degintersection}, as in the previous lemmas, 
\[
\Phi_{n, w}^*\Theta\bullet\Theta=\deg\kappa_4(\Psi_{n, w}(x,y))=\max\{\deg \mu_{1,n}^{ w},\deg\mu_{2,n}^{ w}\},
\]
which proves the lemma for $\Theta_n$ a pole of $\kappa_4$. Similarly if $\Theta_n$ is a zero of $\kappa_4$ one can translate $\Theta_n$ by some Weierstrass point of $\HC$ in order to avoid the zero or the pole divisor of $\kappa_4$. Therefore the degree of $\kappa_4(\Psi_{n}^{ w}(x,y))$ is well defined.
\end{proof}
Using the previous Lemmas, we know that there is always a $(w,0)\in\HC$ such that the value  $\Theta_n\bullet\Theta$ can be obtained using the degree of the rational function $\kappa_4$ restricted to the generic point of $\Theta_n^w$. Further, we know that if $\Psi_n$ is non-constant but $\kappa_4(\Psi_n(x,y))$ is constant then it must be $0$ or $\infty$ and $\Theta_n=\text{Im}\Psi_n\in\text{Supp}\text{ div}(\kappa_4)$ by Lemma \ref{k4constant}. Hence we have a definition of $\delta_n$ for all $n\in \Z$:
\begin{equation}\label{deltan}
    \delta_n := \Bigg\{\begin{array}{lr}
    0 & \text{if } \Psi_n:\HC\to\ja\text{ is constant};\\
        \max\{\deg(\mu_{1,n}^w),\deg(\mu_{2,n}^w)\}& \text{if } \Theta_n\in\text{Supp div}(\kappa_4); \\
        \max\{\deg(\mu_{1,n}),\deg(\mu_{2,n})\} & \text{otherwise.}   
        \end{array}
\end{equation}
Before proving our \textit{basic identity} for genus $2$, namely $\delta_{n-1}+\delta_{n+1}=2\delta_n+4$ to obtain the Hasse-Weil inequality in this case \textit{\`a la} Manin, we recall an additional result. 
\begin{theo.}[Theorem of the cube for Abelian varieties]\label{fundcor}
Let $\ab$ be an Abelian variety, $\alpha,\beta,\gamma\in\text{\rm End}_{k}(\ab)$ and $D\in\text{Div}(\ab)$, then  
\begin{equation*}
(\alpha+\beta+\gamma)^*D-(\alpha+\beta)^*D-(\alpha+\gamma)^*D-(\beta+\gamma)^*D+\alpha^*D+\beta^*D+\gamma^*D\sim 0
\end{equation*}
where $\sim$ denotes linear equivalence.
\end{theo.}
\begin{proof}
See \cite[Corollary 5.3]{milneAV}.
\end{proof}
\begin{cor.}\label{nmapbundle}
Let $\Theta\in\text{\rm Div}(\ja)$ and consider the multiplication-by-$n$ map $[n]\in\text{\rm End}(\ja)$. Then 
$[n]^*\Theta\sim n^2\Theta\in\text{\rm Div}(\ja)$.  
\end{cor.}
\begin{proof}
This is an application of Theorem \ref{fundcor} with $\alpha=[1],\beta=[-1]$ and $\gamma=[n]$ together with induction w.r.t.\ $n$. For details, see \cite[Corollary 5.4]{milneAV} and use that in the present case $[-1]^*\Theta=\Theta$.
\end{proof}

\begin{theo.}\label{recu2}
Let $\HC$ be a hyperelliptic curve of genus $2$ over $\F_q$ with one rational point at infinity, then:
\begin{equation}\label{recurrenceg2}
\delta_{n-1}+\delta_{n+1}=2\delta_{n}+4.
\end{equation}
Moreover, $\delta_n=2n^2 + n(q+1-\#\HC(\F_q))+2q$.
\end{theo.}
\begin{proof}
As before, we denote $\Phi_{m}:=\phi+[m]\in\text{End}(\ja)$ where $\phi$ is the $q$-th Frobenius map.
We begin with some cases when either $\Psi_{n}$ or $\Psi_{n\pm 1}$ is constant.\\

Suppose that $\Psi_n$ is constant, then by definition (\ref{deltan}) we have that $\delta_n=0$. By Corollary \ref{constantmap} we have that $\Theta_{n\pm 1}=\Theta$ and $\Psi_{n\pm 1}=\pm [1]\circ\iota$. It follows by Lemma \ref{selftheta} and the symmetry of $\Theta$ with respect to $[-1]$ that $\delta_{n\pm 1}=[\pm1]^*\Theta\bullet\Theta=\Theta\bullet\Theta=2$ and the lemma follows.\\ 

Now suppose that $\Psi_{n-1}$ is constant. Then $\delta_{n-1}=0$ and by Corollary \ref{constantmap} we have that $\Theta_{n}=\Theta$ and $\Psi_n=\iota$. An analogous argument as in the previous case shows that $\delta_n=2$.\\
To prove that $\delta_{n+1}=8$, note that by Proposition \ref{exceptional} we have that $\phi=-[n-1]\in\text{End}(\ja)$. Therefore $\Phi_{n+1}:=\phi+[n+1]=[2]$ which means that $\Psi_{n+1}=\Phi_{n+1}\circ\iota=[2]\circ\iota$. Hence by Corollary \ref{nmapbundle} and Lemma \ref{degintersection} we obtain $\delta_{n+1}=\Phi_{n+1}^*\Theta\bullet\Theta=[2]^*\Theta\bullet\Theta=4\Theta\bullet\Theta=8$.\\

The case that $\Psi_{n+1}$ is constant is similar to the previous case: one uses the symmetry of $\Theta$ with respect to $[-1]\in\text{End}(\ja)$ to obtain $\delta_{n-1}=[-2]^*\Theta\bullet\Theta=[2]^*\Theta\bullet\Theta=8$.\\

Now we assume that $\Psi_{n\pm 1}$ and $\Psi_n$ are non-constant.
In the case that $\Theta_n\in\text{Supp}\text{ div}(\kappa_4)$, by Lemma \ref{thetanintheta} there is a $(w,0)\in\HC$ such that $\delta_n=\Phi_{n,w}^*\Theta\bullet\Theta=\deg\kappa_4(\Psi_n^w(x,y))$.  Using Remark \ref{remtrans} and Lemma \ref{degintersection}, the latter integer equals $\Phi_{n}^*\Theta\bullet\Theta$. Note that also in the case
that $\Theta_n\not\in\text{Supp}\text{ div}(\kappa_4)$ we have
$\delta_n=\Phi_n^*\Theta\bullet\Theta$
(see Lemma~\ref{degintersection}).
We will finish the proof by studying these
intersection numbers.\\

Using Theorem \ref{fundcor}, let $D:=\Theta\in\text{Div}(\JacH)$ and take $\alpha:=\phi+[n], \beta=[1],\gamma:=-[1]\in\text{\rm End}_{\F_q}(\JacH)$. 
The Theorem of the cube (\ref{fundcor}) implies:
\begin{equation*}
2\Phi_{n}^*\Theta- \Phi_{n+1}^*\Theta-\Phi_{n-1}^*\Theta + 2\Theta\sim 0,
\end{equation*}
or equivalently:
\begin{equation}\label{recthetapull}
2 \Phi_{n}^*\Theta + 2\Theta\sim \Phi_{n-1}^*\Theta + \Phi_{n+1}^*\Theta.
\end{equation}
Intersecting both sides of the equivalence with $\Theta$ proves the first part of the theorem.  To be more precise, we use Lemma \ref{degintersection} together with Lemma \ref{selftheta} to deduce $2\delta_{n}+4=\delta_{n-1}+\delta_{n+1}$.
\\ 

The explicit formula for $\delta_n$ now
follows by induction, noting
(Proposition~\ref{corok41}) that $\delta_{1}=3(q+1)-\#\HC(\F_q)$
and (Lemma~\ref{selfthetafrob}) that $\delta_0=2q$.
\end{proof}

\begin{cor.}[Hasse-Weil for $g\!=\!2$]
Let $\HC/\F_q$ be a hyperelliptic curve with one rational point at infinity and 
$\mbox{char}(\F_q)\neq 2$, then:
\begin{equation}
|q+1-\#\HC({\F_q})|\leq 4\sqrt{q}.
\end{equation}
\end{cor.}
\begin{proof}
Consider the polynomial in $n$ appearing in the previous Theorem \ref{recu2}. The polynomial has the form $\delta(x):=2x^2+Tx+2q$ with
$T:=q+1-\#\HC(\F_q)$. Its discriminant is 
\[ 
\Delta_\delta:=T^2-16q.
\]
We want to prove that $\Delta_\delta\leq 0$ since that would imply that $|T|\leq 4\sqrt{q}$, which is exactly the statement of the Hasse-Weil inequality for $g=2$.\\ \\ 
We already proved in Proposition \ref{exceptional} that if $n\in\Z$ exists such that $\Psi_n\in\text{Mor}_{\F_q}(\HC,\JacH)$ is constant, then $\Psi_n=0$, $q=n^2$ is a perfect square and $\#\HC(\F_q)=q+1\pm 4\sqrt{q}$. Hence 
from the existence of such $n$, the Hasse-Weil inequality over $\F_q$ for the curve in question follows. So
from now on we will suppose that $\Psi_n$ is non-constant
for every $n\in\Z$. By Theorem~\ref{recu2} this implies that $\delta_n=\Phi_n^*\Theta\bullet\Theta$ for all $n\in\Z$.\\ \\
 
It is clear that $\delta_n>0$ for all $n\in\Z$ by definition. This is since $\Psi_n$ is non-constant, hence $\Theta_n\subset\ja$ is a curve, implying that $\delta_n$ is the degree of the rational function $\kappa_4(\Psi_{n,w}(x,y))\in\F_q(x)^*$ or $\kappa_4(\Psi_{n}(x,y))\in\F_q(x)^*$ depending on $\Theta_n$ being in $\text{Supp div}(\kappa_4)$ or not.

Another fast and not very elementary argument for this uses
that the divisor $\Theta$ is ample, hence by the
Nakai-Moishezon criterion for ampleness on surfaces (see \cite[Chapter V, Theorem 1.10]{HARTSHORNE}), its intersection number with any curve is positive.\\

Now from Theorem \ref{recu2} we have that $\delta_n=2n^2 + (q+1-\#\HC(\F_q))n+2q$. Consider $\delta(x)=2x^2 + (q+1-\#\HC(\F_q))x + 2q$. We claim that $\delta(x)$ is non-negative for all $x\in\mathbb{R}$, hence it has non-positive discriminant $\Delta_\delta$. This will imply the Hasse-Weil inequality for this case.\\

Suppose that the Hasse-Weil inequality for genus $2$ is false. This is equivalent to the statement $\Delta_\delta>0$. In this case $\delta(x)$ has two different real zeros $\alpha<\beta$. We have that $\Delta_\delta$ in terms of $\alpha$ and $\beta$ is given by:
\[
\Delta_\delta=4(\alpha-\beta)^2=T^2-16q.
\]
The integer $\Delta_\delta$ is assumed to be positive,
so we conclude $4(\alpha-\beta)^2\geq 1$.
Moreover, recall $\delta(n)>0$ for every $n\in\mathbb{Z}$. Since for any $x_0\in (\alpha,\beta)$ we have that  $\delta(x_0)<0$, it follows that $(\alpha,\beta)$ contains no integers. This implies that $\beta-\alpha < 1$ and then $1\leq 4(\alpha-\beta)^2< 4$.\\
So we have just three situations for positive discriminant: $T^2-16q\in\{1,2,3\}$. Each of these
possibilities results in a contradiction as we will see below.\\ \\
\textbf{Case $T^2-16q=3$:}
There are no integers $(T,q)$ satisfying this, as one checks by reducing modulo $4$.\\ \\
\textbf{Case $T^2-16q=2$:}
Again reducing modulo $4$ implies that
no integral solutions $(T,q)$ exist.\\ \\
\textbf{Case $T^2-16q=1$:}
Then 

\textbf{Subcase (i): $T=8w+1=q+1-\#\HC(\F_q)$, $q=4w^2+w=p^n$}.

Since $p$ is the only prime dividing $q=w(4w+1)$ and since $\text{gcd}(w,4w+1)=1$, it follows that $w=\pm 1$ or $4w+1=\pm 1$. We proceed to check all possibilities.\\ \\If $4w+1=+1$ then $w=0$ and $q=0$ which is not possible.\\
If $4w+1=-1$ then $w=-\tfrac{1}{2}$  which is absurd
since $w$ is an integer.\\
If $w=+1$ then $q=5$ and $T=9$. 
However $9=5+1-\#\HC(\F_5)$ is impossible.\\
If $w=-1$ then $q=3$ and $T=-7$.
However $\HC(\F_3)$ has at
most $2\cdot 3+2$ rational points, hence
$T\geq 3+1-8=-4$.\\

\textbf{Subcase (ii): $T=8w+7=q+1-\#\HC(\F_q)$, $q=4w^2+7w+3=p^n$}.

Again $p$ is the only prime dividing $q=4w^2+7w+3=(w+1)(4w+3)$. Moreover these two factors are coprime since $4(w+1)-(4w+3)=1$. Therefore one of the factors must be $\pm 1$. Again we check all possibilities.\\

If $w+1=1$ then $q=3$ and $T=7$. However any curve
$C/\F_3$ has at least $0$ rational points, hence
$T=3+1-\#C(\F_3)\leq 4$.\\
If $w+1=-1$ then $q=5$ and $T=-9$.
Any hyperelliptic $\HC/\F_5$ satisfies $\#\HC(\F_5)\leq 2\cdot(5+1)$, hence $T\geq 6-12=-6$.\\
The case $4w+3=1$ is impossible since $w$ is assumed to be an integer.\\
Finally, $4w+3=-1$ leads to $q=0$ which is absurd.
 \\ \\
This shows that assuming $\Delta_\delta=T^2-16q> 0$
leads to a contradiction. Therefore $|T|\leq 4\sqrt{q}$ which is the Hasse-Weil inequality for this case. 
\end{proof}

\noindent {\bf Remark}. Note that our {\em definition} of the integers
$\delta_n$ is elementary and completely analogous to the definition by Manin of
the integers $d_n$. However, whereas Manin also succeeded in presenting
a completely elementary proof of the basic identity for the $d_n$, we used
the interpretation of the $\delta_n$ as intersection numbers in order to
show an analogous basic identity in the genus two case. To obtain a fully
elementary proof also in genus two, it therefore suffices to replace this
intersection theory argument by a calculation in the spirit of what Manin did.
We do not know whether this is a feasible task.
    
\bibliographystyle{apa}

\begin{thebibliography}{10}

\bibitem{PROLEGOMENA}
John William~Scott Cassels and E~Victor Flynn.
\newblock {\em Prolegomena to a middlebrow arithmetic of curves of genus 2},
  volume 230.
\newblock Cambridge University Press, 1996.

\bibitem{chahal1988equations}
J.~S. Chahal.
\newblock Chapter~8: Equations over finite fields.
\newblock In {\em Topics in Number Theory}, pages 147--162. Plenum Press, New
  York, 1988.

\bibitem{chahal2014supplement}
Jasbir~S Chahal, Afzal Soomro, and Jaap Top.
\newblock A supplement to {M}anin's proof of the {H}asse inequality.
\newblock {\em Rocky Mountain Journal of Mathematics}, 44(5):1457--1470, 2014.

\bibitem{FLYNNJAC}
Eugene~Victor Flynn.
\newblock The {J}acobian and formal group of a curve of genus 2 over an
  arbitrary ground field.
\newblock 107(03):425--441, 1990.

\bibitem{intfulton}
William Fulton.
\newblock {\em Intersection theory}, volume~2 of {\em Ergebnisse der Mathematik
  und ihrer Grenzgebiete}.
\newblock Springer-Verlag, Berlin, 1984.

\bibitem{HARTSHORNE}
R.~Hartshorne.
\newblock {\em Algebraic Geometry}, volume~52 of {\em Graduate Texts in
  Mathematics}.
\newblock Springer Verlag, New York, 1977.

\bibitem{Hasse1933}
H.~Hasse.
\newblock Beweis des analogons der riemannschen vermutung für die artinschen
  und f. k. schmidtschen kongruenzzetafunktionen in gewissen elliptischen
  fällen. vorläufige mitteilung.
\newblock {\em Nachrichten von der Gesellschaft der Wissenschaften zu
  Göttingen, Mathematisch-Physikalische Klasse}, 1933:253--262, 1933.

\bibitem{Hasse1934}
Helmut Hasse.
\newblock Abstrakte begr{\"u}ndung der komplexen multiplikation und riemannsche
  vermutung in funktionenk{\"o}rpern.
\newblock {\em Abhandlungen aus dem Mathematischen Seminar der Universit{\"a}t
  Hamburg}, 10(1):325--348, Jun 1934.

\bibitem{hasse1936}
Helmut Hasse.
\newblock Zur {T}heorie der abstrakten elliptischen {F}unktionenk\"orper
  {I,II,III}. {D}ie {S}truktur der {G}ruppe der {D}ivisorenklassen endlicher
  {O}rdnung.
\newblock {\em J. Reine Angew. Math.}, 175:55--62, 1936.

\bibitem{Lang}
Serge Lang.
\newblock {\em Abelian Varieties}, volume~7 of {\em Interscience Tracts in Pure
  and Applied Mathematics}.
\newblock Interscience Publishers, New York, 1959.

\bibitem{manin1956oncubic}
Yu.~I. Manin.
\newblock On cubic congruences to a prime modulus.
\newblock {\em Izv. Akad. Nauk SSSR. Ser. Mat.}, 20:673--678, 1956.

\bibitem{manin1960oncubic}
Yu.~I. Manin.
\newblock On cubic congruences to a prime modulus.
\newblock {\em Amer. Math. Soc. Transl. (2)}, 13:1--7, 1960.

\bibitem{milneAV}
James~S. Milne.
\newblock Abelian varieties (v2.00), 2008.
\newblock {\tt http://www.jmilne.org/math/}.

\bibitem{moretcurves}
Laurent Moret-Bailly.
\newblock Familles de courbes et de vari{\'e}t{\'e}s ab{\'e}liennes sur {P}1.
\newblock {\em S{\'e}m. sur les pinceaux de courbes de genre au moins deux (ed.
  L.~Szpiro). Ast{\'e}risques}, 86, 1981.

\bibitem{Soeten}
M.M.J Soeten.
\newblock {\em Hasse$'$s Theorem on Elliptic Curves}.
\newblock Rijks\-universiteit Groningen, 2013.
\newblock {\tt \url{http://fse.studenttheses.ub.rug.nl/id/eprint/10999}}
  (Master's thesis).

\bibitem{ACOFF}
Muhammad~Afzal Soomro.
\newblock {\em Algebraic curves over finite fields}.
\newblock Rijksuniversiteit Groningen, 2013.
\newblock {\tt
  \url{http://hdl.handle.net/11370/024430b9-3e8e-497f-8374-326f014a26e7}} (PhD
  Thesis).

\bibitem{wash}
Lawrence~C. Washington.
\newblock {\em Elliptic curves Number theory and cryptography}.
\newblock Discrete Mathematics and its Applications. Chapman \& Hall / CRC,
  Boca Raton, second edition, 2008.

\end{thebibliography}

\end{document}